\documentstyle{article}

\input{tcilatex}
\begin{document}

\title{{\bf Abstraction and Application in Adjunction}}
\author{{\sc Kosta Do\v sen} \\
Matemati\v {c}ki institut, SANU\\
Knez Mihailova 35, p.f. 367 \\
11001 Belgrade, Yugoslavia \\
email: kosta@mi.sanu.ac.yu}
\date{}
\maketitle

\begin{abstract}
\noindent The postulates of comprehension and extensionality in set theory
are based on an inversion principle connecting set-theoretic abstraction and
the property of having a member. An exactly analogous inversion principle
connects functional abstraction and application to an argument in the
postulates of the lambda calculus. Such an inversion principle arises also
in two adjoint situations involving a cartesian closed category and its
polynomial extension. Composing these two adjunctions, which stem from the
deduction theorem of logic, produces the adjunction connecting product and
exponentiation, i.e. conjunction and implication.\\[0.2cm]
\end{abstract}

\section{Introduction}

\noindent If one bases set theory on two notions, one being abstraction of a
set from a property, i.e. finding the extension of the property, and the
other the property of having a member, then the fundamental postulates of
comprehension and extensionality may be understood as stating that these two
notions are inverse to each other. These two set-theoretical postulates are
analogous to the postulates of beta and eta conversion in the lambda
calculus, where the role of set abstraction is played by functional
abstraction, and the role of having a member by application to an argument.
Abstraction binds a variable and application to a variable introduces it.

An analogous inversion principle arises also in two adjoint situations
involving a cartesian closed category and its polynomial extension. In one
of these adjunctions we find for the functor that maps the original
cartesian closed category to its image in the polynomial extension a
left-adjoint functor based on product, and in the other we find for this
functor a right-adjoint functor based on exponentiation. These two
adjunctions, which stem ultimately from the deduction theorem of logic, and
which had been anticipated in combinatory logic, were first recognized by
Lambek under the name {\it functional completeness} in his pioneering work
in categorial proof theory (see \cite{Lam74}, \cite{LS86}, Part I, and
references therein). Functional completeness is presented quite explicitly
as adjunction in \cite{Jac95} and \cite{D.96}.

After a preliminary section on matters pertaining to the inversion principle
of the postulates of set theory and of the lambda calculus, we shall turn to
categorial proof theory and cartesian closed categories. We shall review the
construction of a polynomial extension of a cartesian closed category,
because this construction, though not difficult, is usually not presented
with sufficient accuracy and detail. Then we shall go through the main steps
of the proof of the two adjunctions of functional completeness, one
involving product and the other exponentiation. We shall see that when these
two adjunctions are composed they give the usual adjunction connecting
product and exponentiation in cartesian closed categories, which is
well-known from Lawvere's work \cite{Law69}.

We shall define precisely all we need, but we shall omit the calculations in
proofs. These calculations are not entirely trivial, but they would not be
new, and a reader with some previous experience with cartesian closed
categories (which he may have acquired by reading, for example, \cite{LS86}%
), or with categories in general (for which many rely on \cite{McL71}),
should be able to perform them.

\section{Set-Theoretical Postulates and Lambda Conversion}

\noindent The two grammatical categories of terms (i.e. individual terms)
and of propositions are basic grammatical categories, with whose help other
grammatical categories can be defined as functional categories: predicates
map terms into propositions, functional expressions map terms into terms,
and connectives and quantifiers map propositions into propositions.

The set-abstracting expression $\{x:...\}$ maps a proposition $A$ into the
term $\{x:A\}$, where the variable $x$ is bound. This term is significant in
particular when $x$ is free in $A$, but it makes sense for any $A$ too. The
expression $x\in ...$ is a unary predicate: it maps a term $a$ into the
proposition $x\in a$. The ideal set theory would just assume that $\{x:...\}$
and $x\in ...$ are inverse to each other, according to the following
postulates: 
\begin{eqnarray*}
&&\text{{\it Comprehension}:\quad }x\in \{x:A\}\leftrightarrow A, \\
&&\text{{\it Extensionality}:\qquad }\{x:x\in a\}=a,
\end{eqnarray*}
provided $x$ is not free in $a$.

In the presence of replacement of equivalents and of Comprehension,
Extensionality is equivalent to the more usual extensionality postulate 
\[
\text{{\it Extensionality}*:\quad }\forall x(x\in a_{1}\leftrightarrow x\in
a_{2})\rightarrow a_{1}=a_{2}, 
\]
provided $x$ is not free in $a_{1}$ and $a_{2}$. The replacement of
equivalents needed here is the principle that from $\forall
x(A_{1}\leftrightarrow A_{2})$ we can infer $\{x:A_{1}\}=\{x:A_{2}\}$, which
can be understood as a principle of logic. That Extensionality entails
Extensionality* is shown as follows. From the antecedent of Extensionality*
with replacement of equivalents we obtain $\{x:x\in a_{1}\}=\{x:x\in a_{2}\}$%
, which yields $a_{1}=a_{2}$ by Extensionality. To show that, conversely,
Extensionality* entails Extensionality, we have $x\in \{x:x\in
a\}\leftrightarrow x\in a$ by Comprehension, from which we obtain
Extensionality by universal generalization and Extensionality*. If doubt is
cast on the replacement of equivalents used above, note that this principle
is implied by Comprehension and Extensionality*. From $\forall
x(A_{1}\leftrightarrow A_{2})$ by Comprehension we obtain $\forall x(x\in
\{x:A_{1}\}\leftrightarrow x\in \{x:A_{2}\})$, and then by Extensionality*
we obtain $\{x:A_{1}\}=\{x:A_{2}\}$.

With the help of substitution for free variables, which is also a principle
of pure logic, we derive the following form of Comprehension: 
\[
\text{{\it Comprehension}*:\quad }y\in \{x:A\}\leftrightarrow A_{y}^{x}, 
\]
where $A_{y}^{x}$ is obtained by substituting uniformly $y$ for free
occurrences of $x$ in $A$, provided the usual provisos for substitution are
satisfied. These provisos will be satisfied if $y$ doesn't occur in the
proposition $A$ at all, neither free nor bound. For such a $y$ we have by
Extensionality 
\[
\{x:A\}=\{y:y\in \{x:A\}\}, 
\]
which with Comprehension* and replacement of equivalents gives $%
\{x:A\}=\{y:A_{y}^{x}\}.$

We know, of course, that ideal set theory is inconsistent if in propositions
we find negation, or at least implication. To get consistency, either $%
\{x:A\}$ will not always be defined, and we replace Comprehension by a
number of restricted postulates, or we introduce types for terms.

Instead of $\{x:...\}$ let us now write $(\lambda _{x}...)$, and instead of $%
x\in ...$ let us write $(...x)$. Then Comprehension and Extensionality
become respectively 
\begin{eqnarray*}
((\lambda _{x}A)x) &\leftrightarrow &A, \\
(\lambda _{x}(ax)) &=&a.
\end{eqnarray*}

If we take that $(\lambda _{x}...)$ maps a term $a$ into the term $(\lambda
_{x}a)$, while $(...x)$ maps a term $a$ into the term $(ax)$, and if,
furthermore, we replace equivalence by equality, and omit outermost
parentheses, our two postulates become the following postulates of the
lambda calculus: 
\begin{eqnarray*}
&&\beta \text{{\it -equality}:\quad }(\lambda _{x}a)x=a, \\
&&\eta \text{{\it -equality}:\quad }\lambda _{x}(ax)=a,
\end{eqnarray*}
provided $x$ is not free in $a$ in $\eta $-equality. The present form of $%
\beta $-equality yields the usual form 
\[
(\lambda _{x}a)b=a_{b}^{x} 
\]
in the presence of substitution for free variables. The usual form of $\beta 
$-equality and $\eta $-equality imply the $\alpha $-equality $\lambda
_{x}a=\lambda _{y}a_{y}^{x}$, provided $y$ doesn't occur in $a$; we proceed
as in the derivation of $\{x:A\}=\{y:A_{y}^{x}\}$ above. The fact that the
lambda calculus based on $\beta $-equality and $\eta $-equality is
consistent is due to the fact that the language has been restricted, either
by preventing anything like negation or implication to occur in terms, or by
introducing types. Without restrictions, in type-free illative theories, we
regain inconsistency.

So the general pattern of Comprehension and Extensionality, on the one hand,
and of $\beta $ and $\eta $-equality, on the other, is remarkably analogous.
These postulates assert that a variable-binding, abstracting, expression $%
\Gamma _{x}$ and application to a variable $\Phi _{x}$ are inverse to each
other, in the sense that $\Phi _{x}\Gamma _{x}\alpha $ and $\Gamma _{x}\Phi
_{x}\alpha $ are either equivalent or equal to $\alpha $, depending on the
grammatical category of $\alpha $. It is even more remarkable that theories
so rich and important as set theory and the lambda calculus are based on
such a simple inversion principle.

\section{The Deduction Theorem in Categorial Proof Theory}

\noindent To speak about deductions we may use labelled sequents of the form 
$f:\Gamma \vdash B$, where $\Gamma $ is a collection of propositions making
the premises, the proposition $B$ is the conclusion, and the term $f$
records the rules justifying the deduction. If the premises can be collected
into a single proposition, and this is indeed the case if $\Gamma $ is
finite and we have a connective like conjunction, then we can restrict our
attention to simple sequents of the form $f:A\vdash B$, where both $A$ and $%
B $ are propositions. We can take that $f:A\vdash B$ is an arrow in a
category in which $A$ and $B$ are objects.

Special arrows in a category are axioms, and operations on arrows are rules
of inference. Equalities of arrows are equalities of deductions. For that,
categorial equalities between arrows have to make proof-theoretical sense,
as indeed they do in many sorts of categories, where they follow closely
reductions in a normalization or cut-elimination procedure. In particular,
equalities between arrows in cartesian closed categories correspond to
equivalence between deductions induced by normalization or cut-elimination
in the implication-conjunction fragment of intuitionistic logic.

Our purpose here is to show that in the context of deductions, as they are
understood in categories, there is something analogous to the inversion
principle we encountered before in set theory and the lambda calculus.

Take a category ${\cal K}$ with a terminal object {\sf T} (this object
behaves like the constant true proposition), and take the polynomial
category ${\cal K}[x]$ obtained by extending ${\cal K}$ with an
indeterminate arrow $x:{\sf T}\vdash D$. Below, we shall explain precisely
what this means, but let us introduce this matter in a preliminary manner.
We obtain ${\cal K}[x]$ by adding to the graph of arrows of ${\cal K}$ a new
arrow $x:{\sf T}\vdash D$, and then by imposing on the new graph equalities
required by the particular sort of category to which ${\cal K}$ belongs.
Note that ${\cal K}[x]$ is not simply the free category of the required sort
generated by the new graph, because the operations on objects and arrows of $%
{\cal K}[x]$ should coincide with those of ${\cal K}$ on the objects and
arrows inherited from ${\cal K}$. We can conceive of ${\cal K}[x]$ as the
extension of a deductive system ${\cal K}$ with a new axiom $D$.

Now consider the variable-binding expression $\Gamma _{x}$ that assigns to
every arrow term $f:A\vdash B$ of ${\cal K}[x]$ the arrow term $\Gamma
_{x}f:A\vdash D\rightarrow B$ of ${\cal K}$, where $\rightarrow $, which
corresponds to implication, is a binary total operation on the objects of $%
{\cal K}$ (in categories, $D\rightarrow B$ is more often written $B^{D}$).
Passing from $f$ to $\Gamma _{x}f$ corresponds to the deduction theorem.
Conversely, we have application to $x$, denoted by $\Phi _{x}$, which
assigns to an arrow term $g:A\vdash D\rightarrow B$ of ${\cal K}$ the arrow
term $\Phi _{x}g:A\vdash B$ of ${\cal K}[x]$. Now, passing from $g$ to $\Phi
_{x}g$ corresponds to modus ponens.

If we require that 
\begin{eqnarray*}
(\beta )\qquad \Phi _{x}\Gamma _{x}f &=&f, \\
(\eta )\qquad \Gamma _{x}\Phi _{x}g &=&g,
\end{eqnarray*}
we obtain a bijection between the hom-sets ${\cal K}[x](A,B)$ and ${\cal K}%
(A,D\rightarrow B)$. If, moreover, we require that this bijection be natural
in the arguments $A$ and $B$, we obtain an adjunction. The left-adjoint
functor in this adjunction is the {\it heritage} functor from ${\cal K}$ to $%
{\cal K}[x]$, which assigns to objects and arrows of ${\cal K}$ their heirs
in ${\cal K}[x]$, while the right-adjoint functor is a functor from ${\cal K}%
[x]$ to ${\cal K}$ that assigns to an object $B$ the object $D\rightarrow B$%
. We find such an adjunction in cartesian closed categories, whose arrows
correspond to deductions of the implication-conjunction fragment of
intuitionistic logic, and also in bicartesian closed categories, whose
arrows correspond to deductions of the whole of intuitionistic propositional
logic. (In bicartesian categories we have besides all finite products,
including the empty product, i.e. terminal object, all finite coproducts,
including the empty coproduct, i.e. initial object.)

In cartesian closed and bicartesian closed categories, as well as in
cartesian categories tout court (namely, in categories with all finite
products), we also have the adjunction given by the bijection between the
hom-sets ${\cal K}(D\times A,B)$ and ${\cal K}[x](A,B)$. Here the heritage
functor is right adjoint, and a functor from ${\cal K}[x]$ to ${\cal K}$
that assigns to an object $A$ the object $D\times A$ is left adjoint. The
binary product operation on objects $\times $\ corresponds to conjunction,
both intuitionistic and classical, as $\rightarrow $\ corresponds to
intuitionistic implication.

Actually, in cartesian closed categories we don't need the terminal object
to express the adjunction involving $\rightarrow $. We could as well take an
indeterminate $x:C\vdash D$, and show that there is a bijection between the
hom-sets ${\cal K}[x](A,B)$ and ${\cal K}(A,(C\rightarrow D)\rightarrow B)$,
natural in the arguments $A$ and $B$. Such an adjunction could also be
demonstrated for categories that have only exponentiation and lack product.
These categories, which correspond to the lambda calculus with only
functional types, are not usually considered. This is probably because their
axiomatization is not very transparent. It is similar to axiomatizations of
systems of combinators \`{a} la Sch\"{o}nfinkel and Curry, where to catch
extensionality we have some rather unwieldy equalities. The main difference
with the axiomatizations of systems of combinators is that in categories
composition replaces functional application, but otherwise these
axiomatizations are analogous.

These adjunctions, which are a refinement of the deduction theorem, were
first considered by Lambek under the name {\it functional completeness} (see
references above; in his first paper on functional completeness \cite{Lam72}
Lambek actually envisaged rather unwieldy combinatorially inspired
equalities, like those we mentioned in the previous paragraph). Through the
categorial equivalence of the typed lambda calculus with cartesian closed
categories, which was discovered by Lambek in the same papers, our
adjunctions are closely related to the so-called Curry-Howard correspondence
between typed lambda terms and natural-deduction proofs. They shed much
light on this correspondence.

\section{Cartesian Closed Categories}

\noindent Although it is assumed the reader has already some acquaintance
with categories, and with cartesian closed categories in particular, to fix
notation and terminology we have to go through some elementary definitions.

A {\it graph} is a pair of functions, called the {\it source} and {\it target%
} function, from a set whose members are called {\it arrows} to a set whose
members are called {\it objects}. (We speak only of small graphs, and small
categories later.) We use $f,g,h,...$, possibly with indices, for arrows,
and $A,B,C,...$, possibly with indices, for objects. We write $f:A\vdash B$
to say that $A$ is the source of $f$ and $B$ its target; $A\vdash B$ is the 
{\it type} of $f$. (We write the turnstile $\vdash $ instead of the more
usual $\rightarrow $, which we use below instead of exponentiation.)

A {\it deductive system} is a graph in which for every object $A$ we have a
special arrow ${\bf 1}_{A}$ $:A\vdash A$, called an {\it identity arrow},
and whose arrows are closed under the binary partial operation of {\it %
composition}: 
\[
\frac{f:A\vdash B\qquad \qquad g:B\vdash C}{g\circ f:A\vdash B} 
\]

A {\it category} is a deductive system in which the following {\it %
categorial equalities} between arrows are satisfied: 
\begin{eqnarray*}
&&f\circ {\bf 1}_{A}={\bf 1}_{B}\circ f=f, \\
&&h\circ (g\circ f)=(h\circ g)\circ f.
\end{eqnarray*}

A {\it cartesian closed deductive system}, or {\it CC} {\it system} for
short, is a deductive system in which we have a special object ${\sf T}$,
and the objects are closed under the binary total operations on arrows $%
\times $ and $\rightarrow $; moreover, for all objects $A$, $A_{1}$, $A_{2}$%
, $B$ and $C$ we have the special arrows (i.e. nullary operations) 
\begin{eqnarray*}
&&k_{A}:A\vdash {\sf T}, \\
&&p_{A_{1},A_{2}}^{i}:A_{1}\times A_{2}\vdash A_{i},~~for~~i\in \{1,2\}, \\
&&\varepsilon _{A,B}:A\times (A\rightarrow B)\vdash B,
\end{eqnarray*}
and the partial operations on arrows 
\[
\frac{f_{1}:C\vdash A_{1}\qquad \qquad f_{2}:C\vdash A_{2}}{\langle
f_{1},~f_{2}\rangle :C\vdash A_{1}\times A_{2}} 
\]
\[
\frac{f:A\times C\vdash B}{\gamma _{A,C}f:C\vdash A\rightarrow B} 
\]

We shall find it handy to use the following abbreviations:

\begin{eqnarray*}
&&\text{{\it for}}{\it ~}f:A\vdash B{\it ~}\text{{\it and}}{\it ~}g:C\vdash
D, \\
&&\quad f\times g=_{def}\langle f\circ p_{A,C}^{1},~g\circ
p_{A,C}^{2}\rangle :A\times C\vdash B\times D, \\
&&\quad f\rightarrow g=_{def}\gamma _{A,B\rightarrow C}(g\circ \varepsilon
_{B,C}\circ (f\times {\bf 1}_{B\rightarrow C})):B\rightarrow C\vdash
A\rightarrow D,\bigskip \\
&&\text{{\it for}}{\it ~}g:C\vdash A\rightarrow B,\quad \varphi
_{A,B}g=_{def}\varepsilon _{A,B}\circ ({\bf 1}_{A}\times g):A\times C\vdash
B,\bigskip \\
&&\quad \overleftarrow{b}_{A,B,C}=_{def}\langle p_{A,B}^{1}\circ p_{A\times
B,C}^{1},~p_{A,B}^{2}\times {\bf 1}_{C}\rangle :(A\times B)\times C\vdash
A\times (B\times C), \\
&&\quad \overrightarrow{b}_{A,B,C}=_{def}\langle {\bf 1}_{A}\times
p_{B,C}^{1},~p_{B,C}^{2}\circ p_{A,B\times C}^{2}\rangle :A\times (B\times
C)\vdash (A\times B)\times C, \\
&&\quad c_{A,B}=_{def}\langle p_{A,B}^{2},~p_{A,B}^{1}\rangle :A\times
B\vdash B\times A.
\end{eqnarray*}

A {\it cartesian closed category}, or {\it CC category} for short, is a CC
system in which besides the categorial equalities the following {\it CC
equalities} hold: 
\begin{eqnarray*}
&&\text{(}{\sf T}\eta \text{)}\qquad \text{{\it for}}{\it \ }f:A\vdash {\sf T%
},\quad f=k_{A},\bigskip \\
&&\text{(}\times \beta \text{)}\,\qquad p_{A_{1},A_{2}}^{i}\circ \langle
f_{1},~f_{2}\rangle =f_{i}, \\
&&\text{(}\times \eta \text{)}\qquad \langle p_{A,B}^{1}\circ
h,~p_{A,B}^{2}\circ h\rangle =h,\bigskip \\
&&\text{(}\rightarrow \beta \text{)}\qquad \varphi _{A,B}\gamma _{A,C}f=f, \\
&&\text{(}\rightarrow \eta \text{)}\qquad \gamma _{A,C}\varphi _{A,B}g=g.
\end{eqnarray*}

If ${\cal K}$ and ${\cal L}$ are CC categories, a {\it strict cartesian
closed functor}, or for short {\it CC functor}, F from ${\cal K}$ to ${\cal L%
}$ is a functor that satisfies the following equalities in ${\cal L}$: 
\begin{eqnarray*}
&&\text{F{\sf T}}=\text{{\sf T}},\quad \text{F}(A~\alpha ~B)=\text{F}%
A~\alpha ~\text{F}B\text{,{\it \ where }}\alpha \text{{\it \ is }}\times 
\text{{\it or} }\rightarrow \text{,} \\
&&\text{F}k_{A}=k_{\text{F}A},\quad \text{F}p_{A,B}^{i}=p_{\text{F}A,\text{F}%
B}^{i},\quad \text{F}\varepsilon _{A,B}=\varepsilon _{\text{F}A,\text{F}B},
\\
&&\text{F}\langle f_{1},~f_{2}\rangle =\langle \text{F}f_{1},~\text{F}%
f_{2}\rangle ,\quad \text{F}\gamma _{A,C}f=\gamma _{\text{F}A,\text{F}C}%
\text{F}f.
\end{eqnarray*}

\section{The Polynomial Cartesian Closed Category}

\noindent Given a CC category ${\cal K}$, and an object $D$ of ${\cal K}$,
we shall construct the {\it polynomial} CC category ${\cal K}[x]$ obtained
by adjoining an {\it indeterminate} arrow $x:{\sf T}\vdash D$ by first
constructing a CC system ${\cal S}$ obtained by adjoining the indeterminate
arrow $x$ to ${\cal K}$.

The objects of ${\cal S}$ will be the same as the objects of ${\cal K}$. We
provide a mathematical object $x$, which is not an arrow of ${\cal K}$, and
different mathematical objects denoted by $\circ ^{{\cal S}}$, $\langle
,\rangle ^{{\cal S}}$ and $\gamma _{A,B}^{{\cal S}}$, for every pair $(A,B)$
of objects of ${\cal K}$. Then we define inductively the {\it arrows of} $%
{\cal S}$:\medskip

\noindent (0) $x$ is an arrow of ${\cal S}$ of type ${\sf T}\vdash D$%
;\smallskip

\noindent (1) every arrow of ${\cal K}$ is an arrow of ${\cal S}$, with the
same type it has in ${\cal K}$;\smallskip

\noindent (2) if $f:A\vdash B$ and $g:B\vdash C$ are arrows of ${\cal S}$,
then the ordered triple $(\circ ^{{\cal S}},f,g)$ is an arrow of ${\cal S}$
of type $A\vdash C$;\smallskip

\noindent (3) if $f_{1}:C\vdash A_{1}$ and $f_{2}:C\vdash A_{2}$ are arrows
of ${\cal S}$, then the ordered triple $(\langle ,\rangle ^{{\cal S}%
},f_{1},f_{2})$ is an arrow of ${\cal S}$ of type $C\vdash A_{1}\times A_{2}$%
;\smallskip

\noindent (4) if $f:A\times C\vdash B$ is an arrow of ${\cal S}$, then the
ordered pair $(\gamma _{A,C}^{{\cal S}},f)$ is an arrow of ${\cal S}$ of
type $C\vdash A\rightarrow B$.\medskip

\noindent We denote $(\circ ^{{\cal S}},f,g)$, $(\langle ,\rangle ^{{\cal S}%
},f_{1},f_{2})$ and $(\gamma _{A,C}^{{\cal S}},f)$ by $g\circ ^{{\cal S}}f$, 
$\langle f_{1},~f_{2}\rangle ^{{\cal S}}$ and $\gamma _{A,C}^{{\cal S}}f$
respectively.

The CC category ${\cal K}[x]$ will have the same objects as ${\cal K}$ and $%
{\cal S}$, while its arrows will be obtained by factoring the arrows of $%
{\cal S}$ through a suitable equivalence relation. Consider the equivalence
relations $\equiv $ on the arrows of ${\cal S}$ that satisfy the congruence
law 
\[
\text{{\it if} }f_{1}\equiv f_{2}\text{{\it \ and }}g_{1}\equiv g_{2}\text{%
{\it , then} }g_{1}\circ ^{{\cal S}}f_{1}\equiv g_{2}\circ ^{{\cal S}}f_{2} 
\]
(provided the types of the arrows on the two sides of $\equiv $ are equal,
and are such that $g_{1}\circ ^{{\cal S}}f_{1}$ is an arrow of ${\cal S}$).
Moreover, these equivalence relations satisfy analogous congruence laws for $%
\langle ,\rangle ^{{\cal S}}$ and $\gamma _{A,C}^{{\cal S}}$, and they
satisfy basic equivalences obtained from the categorial and CC equalities by
replacing the equality sign $=$ by $\equiv $, and by superscribing ${\cal S}$
on $\circ $, $\rangle $ and $\gamma $. Finally, our equivalence relations
satisfy the following basic equivalences for $f$, $g$, $f_{1}$ and $f_{2}$
arrows of ${\cal K}$ of the appropriate types: 
\begin{eqnarray*}
g\circ ^{{\cal S}}f &\equiv &g\circ f, \\
\langle f_{1},~f_{2}\rangle ^{{\cal S}} &\equiv &\langle f_{1},~f_{2}\rangle
, \\
\gamma _{A,C}^{{\cal S}}f &\equiv &\gamma _{A,C}f,
\end{eqnarray*}
where the operations on the right-hand sides are those of ${\cal K}$. Let us
call equivalence relations that satisfy all that {\it CC equivalence
relations} on the arrows of ${\cal S}$.

It is clear that the intersection of all CC equivalence relations on the
arrows of ${\cal S}$ is again a CC equivalence relation on the arrows of $%
{\cal S}$---the smallest such relation---, which we denote by $\equiv _{\cap
}$. Then for every arrow $f$ of ${\cal S}$ take the equivalence class $[f]$
made of all the arrows $f^{\prime }$ of ${\cal S}$ such that $f\equiv _{\cap
}f^{\prime }$.

The objects of ${\cal K}[x]$ are the objects of ${\cal K}$, and its arrows
are the equivalence classes $[f]$, the type of $[f]$ in ${\cal K}[x]$ being
the same as the type of $f$ in ${\cal S}$ (all arrows in the same
equivalence class have the same type in ${\cal S}$). With the definitions 
\begin{eqnarray*}
&&{\bf 1}_{A}=_{def}[{\bf 1}_{A}], \\
&&\,[g]\circ [f]=_{def}[g\circ ^{{\cal S}}f],
\end{eqnarray*}
and other analogous definitions, it is clear that ${\cal K}[x]$ is a CC
category.

Note that ${\cal K}[x]$ is not the same as the free CC category generated by
the graph of ${\cal K}$ extended with $x$. To pass from this free CC
category to ${\cal K}[x]$ involves further factoring of objects and arrows
through suitable equivalence relations, so as to ensure that the new
operations on objects and arrows coincide with the old operations on the
objects and arrows of ${\cal K}$. However, the extension of ${\cal K}$ to $%
{\cal K}[x]$ is free in a certain sense, which we shall explicate in the
next section.

\section{The Heritage Functor}

\noindent We shall now define a CC functor H from ${\cal K}$ to ${\cal K}[x]$%
, which is called the {\it heritage functor}. On objects H is the identity
function, while on arrows it is defined by 
\[
\text{H}f=_{def}[f]. 
\]

This function on arrows is clearly not onto, because of the arrow $x$ and
other arrows of ${\cal K}[x]$ involving $x$. It is also in general not
one-one. Conditions that ensure that H is one-one on arrows are investigated
in \cite{Lam74}, \cite{LS86} (I.5) and, especially, \cite{Cub98}. A
necessary and sufficient condition, found in this last paper, is that the
object $D$ of $x:{\sf T}\vdash D$ be ``nonempty'', nonemptiness being
expressed in a categorial manner by requiring that the arrow $k_{D}:D\vdash 
{\sf T}$ be epi, i.e. cancellable on the right-hand side of compositions. (A
functor such as H, which is a bijection on objects, is full if and only if
it is onto on arrows, and it is faithful if and only if it is one-one on
arrows.)

It is easy to check that H is a CC functor. For example, we have 
\[
\text{H}(g\circ f)=[g\circ f]=[g\circ ^{{\cal S}}f]=[g]\circ [f]=\text{H}%
g\circ \text{H}f, 
\]
and we proceed analogously in other cases.

The polynomial CC category ${\cal K}[x]$ and the heritage functor H satisfy
the following universal property, which explains in what sense the extension
of ${\cal K}$ to ${\cal K}[x]$ is free:

\begin{quotation}
\noindent {\it For every CC category} ${\cal L}${\it , every CC functor} M 
{\it from} ${\cal K}$ {\it to} ${\cal L}$ {\it and every arrow} $f:{\sf T}%
\vdash $M$D$ {\it of} ${\cal L}$, {\it there is a unique CC functor} N {\it %
from} ${\cal K}[x]$ {\it to} ${\cal L}$ {\it such that} N$x=f$ {\it and} M$%
\,=\,$NH.
\end{quotation}

This property characterizes ${\cal K}[x]$ up to isomorphism. It is analogous
to the universal property one finds in the construction of a polynomial ring 
${\cal K}[x]$ by adding an indeterminate $x$ to a commutative ring ${\cal K}$
(see \cite{McLB79}, IV.4). The analogue of the heritage functor is the {\it %
insertion} homomorphism from ${\cal K}$ to ${\cal K}[x]$ (which, however, is
one-one, whereas the heritage functor need not be faithful). This explains
the epithet {\it polynomial} ascribed to ${\cal K}[x]$.

In general, we encounter the same kind of universal property in connection
with variables. The indeterminate $x$ is in fact a variable, and a variable
is a free element, or a free nullary operation. If to the set of terms $%
{\cal A}$ of an algebra of a certain kind we add a variable $x$ so as to
obtain the set of polynomial terms ${\cal A}[x]$, we shall have the
following universal property involving the heritage (or insertion)
homomorphism $h$ from ${\cal A}$ to ${\cal A}[x]$:

\begin{quotation}
\noindent {\it For every algebra} ${\cal B}$ {\it of the same kind as} $%
{\cal A}${\it , every homomorphism} $m$ {\it from} ${\cal A}$ {\it to} $%
{\cal B}$ {\it and every element} $b$ {\it of} ${\cal B}${\it , there is a
unique homomorphism} $n$ {\it from} ${\cal A}[x]$ {\it to} ${\cal B}$ {\it %
such that} $n(x)=b$ {\it and for every element} $a$ {\it of} ${\cal A}$ {\it %
we have} $m(a)\,=n(h(a))\,$.
\end{quotation}

\section{The Heritage Functor has a Left Adjoint}

\noindent Our aim is now to show that the heritage functor H from the CC
category ${\cal K}$ to the polynomial CC category ${\cal K}[x]$ has a left
adjoint.

The arrows of the CC system ${\cal S}$ were defined inductively, and we
shall first define by induction on the complexity of the arrow $f:A\vdash B$
of ${\cal S}$ a function $\Phi _{x}^{\prime }$ that assigns to $f$ the arrow 
$\Phi _{x}^{\prime }f:D\times A\vdash B$ of ${\cal K}$:\medskip

(0)\quad $\Phi _{x}^{\prime }x=p_{D,{\sf T}}^{1},\smallskip $

(1)\quad $\Phi _{x}^{\prime }f=f\circ p_{D,A}^{2},$\ {\it for} $f$ {\it an
arrow of} ${\cal K}$,$\smallskip $

(2)\quad $\Phi _{x}^{\prime }(g\circ ^{{\cal S}}f)=\Phi _{x}^{\prime }g\circ
\langle p_{D,A}^{1},~\Phi _{x}^{\prime }f\rangle ,\smallskip $

(3)\quad $\Phi _{x}^{\prime }\langle f_{1},~f_{2}\rangle ^{{\cal S}}=\langle
\Phi _{x}^{\prime }f_{1},~\Phi _{x}^{\prime }f_{2}\rangle ,\smallskip $

(4)\quad $\Phi _{x}^{\prime }\gamma _{A,C}^{{\cal S}}f=\gamma _{A,D\times
C}(\Phi _{x}^{\prime }f\circ \overleftarrow{b}_{D,A,C}\circ (c_{A,D}\times 
{\bf 1}_{C})\circ \overrightarrow{b}_{A,D,C}).$\medskip

The equalities $[f]=[g]$ of ${\cal K}[x]$ stem from the equivalences $%
f\equiv _{\cap }g$, which can be derived as in a formal system from the
reflexivity of $\equiv _{\cap }$ and the basic equivalences assumed for CC
equivalence relations with the help of replacement of equivalents. We can
prove the following lemma by induction on the length of the derivation of $%
f\equiv _{\cap }g$.

\begin{lemma}
If $[f]=[g]$ in ${\cal K}[x]$, then $\Phi _{x}^{\prime }f=\Phi _{x}^{\prime
}g$ in ${\cal K}$.
\end{lemma}

If we put 
\[
\Phi _{x}^{\prime }[f]=_{def}\Phi _{x}^{\prime }f, 
\]
Lemma 1 guarantees that this defines indeed a function from the arrows of $%
{\cal K}[x]$ to the arrows of ${\cal K}$.

Then we define a function $\Gamma _{x,A}^{\prime }$ that assigns to an arrow 
$f:D\times A\vdash B$ of ${\cal K}$ the arrow $\Gamma _{x,A}^{\prime
}f:A\vdash B$ of ${\cal K}[x]$: 
\[
\Gamma _{x,A}^{\prime }f=_{def}[f]\circ \langle [x]\circ k_{A},~{\bf 1}%
_{A}\rangle . 
\]
Note that here, contrary to what we had in Sections 2 and 3, the
variable-binding function that corresponds to abstraction has $\Phi $ in its
name, while the function that corresponds to application has $\Gamma $.
Before, it was the other way round. We make this switch to conform to the
notation for adjoint situations of \cite{D.96}, \cite{D.99}, and \cite{D.99a}%
. Conforming to this same notation, in the next section matters will return
to what we had in Sections 2 and 3.

We can verify that $\Phi _{x}^{\prime }$ and $\Gamma _{x,A}^{\prime }$
establish a bijection between the hom-sets ${\cal K}(D\times A,B)$ and $%
{\cal K}[x](A,B)$.

\begin{lemma}
For every $[f]:A\vdash B$ of ${\cal K}[x]$ we have $\Gamma _{x,A}^{\prime
}\Phi _{x}^{\prime }[f]=[f]$ in ${\cal K}[x]$.
\end{lemma}

\begin{lemma}
For every $f:D\times A\vdash B$ of ${\cal K}$ we have $\Phi _{x}^{\prime
}\Gamma _{x,A}^{\prime }f=f$ in ${\cal K}$.
\end{lemma}

\noindent We prove Lemma 2 by induction on the complexity of the arrow $f$
of ${\cal S}$ (which involves some not entirely trivial computations when $f$
is of the form $\gamma _{B_{1},D\times A}^{{\cal S}}f^{\prime }$), while
Lemma 3 is checked directly.

We define a functor F from ${\cal K}[x]$ to ${\cal K}$ by 
\begin{eqnarray*}
&&\text{F}A=_{def}D\times A, \\
&&\text{F}[f]=_{def}\Phi _{x}^{\prime }(\Gamma _{x,B}^{\prime }{\bf 1}%
_{D\times B}\circ [f])=\langle p_{D,A}^{1},~\Phi _{x}^{\prime }f\rangle .
\end{eqnarray*}
To check that this is a functor left adjoint to the heritage functor H it
remains to establish 
\[
\Phi _{x}^{\prime }([g]\circ [f])=\Phi _{x}^{\prime }[g]\circ \text{F}[f], 
\]
which was built into the definition of $\Phi _{x}^{\prime }$, and 
\[
\Gamma _{x,A}^{\prime }(f\circ \Phi _{x}^{\prime }{\bf 1}_{A})=[f]=\text{H}f 
\]
(see \cite{D.96}, \S $~$3.1, \cite{D.99}, \S $~$4.1.7, or \cite{D.99a}, \S $%
~ $8).

\section{The Heritage Functor has a Right Adjoint}

\noindent To show that the heritage functor H from the CC category ${\cal K}$
to the polynomial CC category ${\cal K}[x]$ has a right adjoint, we define
first a function $\Gamma _{x}^{\prime \prime }$ that assigns to an arrow $%
[f]:A\vdash B$ of ${\cal K}[x]$ the arrow $\Gamma _{x}^{\prime \prime
}[f]:A\vdash D\rightarrow B$ of ${\cal K}$: 
\[
\Gamma _{x}^{\prime \prime }[f]=_{def}\gamma _{D,A}\Phi _{x}^{\prime }[f]. 
\]

\noindent Then we define a function $\Phi _{x,B}^{\prime \prime }$ that
assigns to an arrow $g:A\vdash D\rightarrow B$ of ${\cal K}$ the arrow $\Phi
_{x,B}^{\prime \prime }g:A\vdash B$ of ${\cal K}[x]$: 
\[
\Phi _{x,B}^{\prime \prime }g=_{def}\Gamma _{x,A}^{\prime }\varphi _{D,B}g. 
\]

\noindent It follows easily from Lemmata 2 and 3, together with the CC
equalities $(\rightarrow ~\beta )$ and $(\rightarrow ~\eta )$, that $\Gamma
_{x}^{\prime \prime }$ and $\Phi _{x,B}^{\prime \prime }$ establish a
bijection between the hom-sets ${\cal K}[x](A,B)$ and ${\cal K}%
(A,D\rightarrow B)$. Namely,

\begin{quotation}
\noindent $(\beta )$\quad {\it for every} $[f]:A\vdash B$ {\it of} ${\cal K}%
[x]$ {\it we have} $\Phi _{x,B}^{\prime \prime }\Gamma _{x}^{\prime \prime
}[f]=[f]$ {\it in} ${\cal K}[x]$;\medskip

\noindent $(\eta )$\quad {\it for every} $g:A\vdash D\rightarrow B$ {\it of} 
${\cal K}$ {\it we have} $\Gamma _{x}^{\prime \prime }\Phi _{x,B}^{\prime
\prime }g=g$ {\it in} ${\cal K}$.
\end{quotation}

We define a functor G from ${\cal K}[x]$ to ${\cal K}$ by 
\begin{eqnarray*}
&&\text{G}A=_{def}D\rightarrow A, \\
&&\text{G}[f]=_{def}\Gamma _{x}^{\prime \prime }([f]\circ \Phi
_{x,A}^{\prime \prime }{\bf 1}_{D\rightarrow A})=\gamma _{D,D\rightarrow
A}(\Phi _{x}^{\prime }f\circ \langle p_{D,D\rightarrow A}^{1},~\varepsilon
_{D,A}\rangle ).
\end{eqnarray*}
To check that this is a functor right adjoint to the heritage functor H it
remains to establish either 
\[
\Gamma _{x}^{\prime \prime }([g]\circ [f])=\text{G}[g]\circ \Gamma
_{x}^{\prime \prime }[f] 
\]
or 
\[
\Phi _{x,C}^{\prime \prime }(g\circ f)=\Phi _{x,C}^{\prime \prime }g\circ
[f], 
\]
together with 
\[
\Phi _{x,B}^{\prime \prime }(\Gamma _{x}^{\prime \prime }{\bf 1}_{B}\circ
f)=[f]=\text{H}f, 
\]
which can be done after some calculation.

Consider now the functors FH and GH from ${\cal K}$ to ${\cal K}$ obtained
by composing the functors F and G, respectively, with the heritage functor
H. It is clear that the functors FH and GH make an adjoint situation in
which FH is left adjoint and GH is right adjoint. This adjunction is the
usual adjunction that ties $D\times $ and $D\rightarrow $ in CC categories,
since we can verify that 
\begin{eqnarray*}
\text{FH}f &=&{\bf 1}_{D}\times f, \\
\text{GH}f &=&{\bf 1}_{D}\rightarrow f.
\end{eqnarray*}
The bijection, natural in the arguments $A$ and $B,$ between ${\cal K}%
(D\times A,B)$ and ${\cal K}(A,D\rightarrow B)$ is given by the operations $%
\gamma _{D,A}$ and $\varphi _{D,B}$. (This bijection is actually natural in
the argument $D$ too.)

\section{Logical Constants and Adjunction}

\noindent Adjointness phenomena pervade logic, as well as much of
mathematics. An essential ingredient of the spirit of logic is to
investigate inductively defined notions, and inductive definitions engender
free structures, which are tied to adjointness. We find also in logic the
important model-theoretical adjointness between syntax and semantics, behind
theorems of the {\it if and only if} type called semantical completeness
theorems. However, adjunction is present in logic most specifically through
its connection with logical constants.

Lawvere put forward the remarkable thesis that all logical constants are
characterized by adjoint functors (see \cite{Law69}). Lawvere's thesis about
logical constants is just one part of what he claimed for adjunction, but it
is a significant part.

For conjunction, i.e. binary product in cartesian categories, we have the
adjunction between the diagonal functor from ${\cal K}$ to the product
category ${\cal K}\times {\cal K}$ as left adjoint and the internal product
bifunctor from ${\cal K}\times {\cal K}$ to ${\cal K}$ as right adjoint.
Properties assumed for this bifunctor are not only sufficient to prove the
adjunction, but they are also necessary---they can be deduced from the
adjunction. Binary coproduct, which corresponds to disjunction, is
analogously characterized as a left adjoint to the diagonal functor. The
terminal and initial objects, which correspond respectively to the constant
true proposition and the constant absurd proposition, may be conceived as
empty product and empty coproduct. They are characterized by functors right
and left-adjoint, respectively, to the constant functor into the trivial
category with a single object and a single identity arrow.

In all that, one of the adjoint functors carries the logical constant to be
characterized, i.e., it involves the corresponding operation on objects, and
depends on the inner constitution of the category, while the other adjoint
functor is a {\it structural} functor, which does not involve the inner
operations of the category, and can be defined for any category
(``structural'' is here used as in the ``structural rules'' of Gentzen's
proof theory). The diagonal functor and the constant functor are clearly
structural: they make sense for any kind of category.

This suggests an amendment to Lawvere's thesis: namely, the functor carrying
the logical constant should be adjoint to a structural functor. This
structural functor is presumably tied to some features of deduction that are
independent of any particular constant we may have in our language, and are
hence {\it formal} in the purest way. With this amendment the thesis might
serve to separate the constants of formal logic from other expressions.

Lawvere's way to characterize intuitionistic implication through adjunction
is by relying on the bijection between ${\cal K}(D\times A,B)$ and ${\cal K}%
(A,D\rightarrow B)$ in cartesian closed categories, which can be obtained by
composing the two adjunctions with the heritage functor, as we have seen in
the previous section. The disadvantage of this characterization is that none
of the adjoint functors $D\times $\ and $D\rightarrow $\ is structural.

Can the adjunctions of functional completeness serve to characterize
conjunction and intuitionistic implication? It would be nice if they could,
because the heritage functor is structural. This is more important for
intuitionistic implication than for conjunction, because for the latter we
already have a characterization through an adjunction with a structural
functor---namely, the adjunction with the diagonal functor. And it would be
preferable if implication were characterized in the absence of conjunction,
and of anything else, as in functional completeness with the categories that
have only exponentiation and may lack product (which we mentioned in Section
3).

To define a polynomial category ${\cal K}[x]$ with an indeterminate $%
x:C\vdash D$ we assume for ${\cal K}[x]$ that it has whatever it must have
to make it a polynomial category of the required kind, to which ${\cal K}$
belongs. This is something, but it is nothing in particular. With an
indeterminate $x:{\sf T}\vdash D$ we assume for the categories ${\cal K}$
and ${\cal K}[x]$ that they have also a special object ${\sf T}$. So, if
intuitionistic implication could be characterized by the adjunction of
functional completeness, this could be achieved even in the absence of ${\sf %
T}$, whereas the characterization of conjunction by the corresponding
adjunction of functional completeness would depend on the presence of the
terminal object ${\sf T}$. We would be able to characterize binary product
only in the presence of empty product, i.e. ${\sf T}$. It seems all finite
products go together. In any case, however, the definition of the heritage
functor is structural. It will be the same for any kind of category.

A step towards showing that conjunction and intuitionistic implication can
be characterized by the adjunctions of functional completeness was taken in 
\cite{D.92}, and, especially, \cite{DP96}. What we need to show is that the
assumptions made for cartesian, or cartesian closed categories, or
categories that have only exponentiation and may lack product, are not only
sufficient for demonstrating the appropriate adjunction of functional
completeness, but they are also necessary. However, here the matter is not
so clear-cut as when the assumptions concerning binary product are deduced
from adjunction with the diagonal functor. It is shown in \cite{DP96}
(Section 5) that many of the assumptions for cartesian categories can be
deduced from functional completeness, but still some assumptions stay simply
postulated. One feels, however, that even these assumptions could be deduced
if matters were formulated in the right way.

\end{document}